\RequirePackage{fix-cm}
\documentclass[smallextended]{svjour3}       % onecolumn (second format)
\smartqed  % flush right qed marks, e.g. at end of proof
\usepackage{graphicx}

\usepackage{subcaption}
\usepackage{comment}
\usepackage{float}
\usepackage{array}
\usepackage{amsmath,amssymb, amsfonts}
\usepackage{url}
\usepackage[font=footnotesize,labelfont=it]{caption}
\usepackage{longtable}
\usepackage{multirow}
\usepackage{xcolor}

\usepackage{enumitem, comment, xifthen}
\newcolumntype{P}[1]{>{\centering\arraybackslash}p{#1}}

\begin{document}
	
	\title{Progressively Strengthening and Tuning MIP Solvers for Reoptimization%\thanks{Grants or other notes
	}
	
	
	\author{Krunal Kishor Patel
	}
	
	\authorrunning{Krunal Patel} % if too long for running head
	
	\institute{K. Patel \at
		CERC, Polytechnique Montr\'eal, Montr\'eal, Quebec, Canada \\
		ORCID 0000-0001-7414-5040\\
		\email{krunal.patel@polymtl.ca}           %  \\
	}
	
	\date{Received: date / Accepted: date}

	\maketitle
	
	\begin{abstract}
		This paper explores reoptimization techniques for solving sequences of similar mixed integer programs (MIPs) more effectively. Traditionally, these MIPs are solved independently, without capitalizing on information from previously solved instances. Our approach focuses on primal bound improvements by reusing the solutions of the previously solved instances, as well as dual bound improvements by reusing the branching history and automating parameter tuning. We also describe ways to improve the solver performance by extending ideas from reliability branching to generate better pseudocosts. {Our reoptimization approach, crafted for the MIP 2023 workshop computational competition, was honored with the first prize.} In this paper, we thoroughly analyze the performance of each technique and their combined impact on the solver's performance. Finally, we present ways to extend our techniques in practice for further improvements.
	\end{abstract}
	
	\section*{Acknowledgements}
	
	The author received financial support from {the} Canada Excellence Research Chair (CERC), Polytechnique Montreal{ and t}he computational resources for experimental evaluation were provided by CERC. The author would also like to thank Prof. Andrea Lodi{,} Prof. Guy Desaulniers{, and the anonymous reviewers} for their helpful comments. 
	
	\section{Introduction}
	
	\begin{comment}

		\begin{itemize}
			\item List of applications where reoptimization can be helpful.
			\item Which components can change.
			\item Past work
			\item Paper organization. (Describe the dual bound high level directions.)
		\end{itemize}
	\end{comment}
	
	{Many applications require solving a sequence of similar mixed integer programs (MIPs). These instances differ on some specific parts (objective, right-hand side, variable bounds). Some of the major examples include decomposition algorithms (that use column generation or row generation) for which the subproblems are similar to each other except changes in the objective or right-hand side (see, e.g., \cite{bolusani2022framework},\cite{witzig2014reoptimization}). Some heuristics (e.g., \cite{danna2005rins}, \cite{berthold2006crossover}, and \cite{gamrath2019structure}) in MIP solvers also solve similar sub-MIPs. This paper describes some techniques to reuse the information generated by solving previous instances to solve future instances more efficiently.}
	
	{In this paper, we use the following formulation to denote a MIP:
	
	\begin{align}
		\begin{split}
			{m}in &\quad c^Tx\\
			s.t. &\quad Ax \leq b\\
			&\quad l \leq x \leq u\\
			&\quad x_j \in \mathbb{Z}, \forall j \in I
		\end{split}
	\end{align}
	
	\noindent
	with $A\in \mathbb{R}^{m \times n}, c \in \mathbb{R}^n, b \in \mathbb{R}^m, l \in (\mathbb{R} \cup \{ -\infty \} )^n, u \in (\mathbb{R} \cup \{ \infty \} )^n$, variables $x \in \mathbb{R}^n$ with $I \subseteq N = \{1,...,n\}$.}
	
	{The research on the topic of reoptimization of MIPs has been limited. Some MIP solvers (e.g.,  Gurobi \cite{gurobi}, CPLEX \cite{cplex2020}) use parts of primal solutions of the previously solved similar instances as warm start. \cite{ausiello2011complexity} and \cite{schieber2018theory} describe reoptimization algorithms for certain classes of combinatorial optimization problems.} \cite{ralphs2006duality} described a way to reuse the dual information. They mainly addressed the changes in the right-hand side. More recently, \cite{gamrath2015reoptimization} investigated techniques for reoptimizaiton and implemented them in the SCIP solver \cite{BestuzhevaEtal2021ZR}, primarily focusing on changes in the objective function. A detailed description can be found in  \cite{witzig2014reoptimization}.
	
	To solve the instances faster, we need to find good quality solutions faster (primal bound improvements) and prove the optimality faster (dual bound improvements). In this paper, we describe primal and dual bound improvement techniques by exploiting the similarity of the MIP instances. We focus on solving some series of 50 similar MIP instances that change in {one or more of the following parts: objective coefficients, constraint bounds, variable bounds, and constraint matrix entries}.
	
	Our reoptimization approach{,} originally developed for the computational competition \cite{bolusani2023mip} of the MIP 2023 workshop{, was honored with} the first prize. The complete source code for this approach is available at \\ https://github.com/krooonal/mipcc2023/tree/MPCPaper. 
	
	{All the techniques described in this paper are inspired by existing literature with minor or no modifications. The major contribution of this paper is using these techniques together to achieve better performance. The techniques include
		
		\begin{enumerate}
			\item Reusing primal solutions: We describe how to combine solutions of the previous instances to form a partial solution hint that the MIP solver can convert into a feasible solution early in the search process. Some of these techniques are inspired by existing heuristics like RINS \cite{danna2005rins}, Crossover \cite{berthold2006crossover}, and CMSA \cite{blum2016cmsa}.
			\item Reusing branching information: MIP solvers progressively learn the effects of branching on each variable in the form of pseudocosts and other statistics, and we reuse these statistics as a warm start for the subsequent instances. This technique is also used in the context of restarts in MIPs \cite{berthold2022restarttransfer}.
			\item Improving pseudocosts: Inspired by the main idea behind reliability branching \cite{achterberg2005reliability}, we explore how using more strong branching in the first instance affects the solution process of the subsequent instances.
			\item Automated parameter tuning: Multi-armed bandit solution techniques are known for various hyperparameter tuning tasks in machine learning (see, e.g., \cite{li2017hyperband}). In MIP solvers, multi-armed bandit algorithms are used in a node selection rule \cite{sabharwal2012guiding}, and adaptive large neighborhood search (ALNS) \cite{hendel2022adaptive}(a primal heuristic for MIP solvers). We show how to use an algorithm for a multi-armed bandit problem (with slight modifications) to tune MIP solver parameters without any offline training. 
			\item Other optimizations: We detect the solver parts that are not working well for a given series of instances and turn them off. This technique can be seen as an extension of parameter tuning.
		\end{enumerate}
	}
	
	This paper focuses on reusing information generated by solving previous instances. We do not explore any offline learning approaches as enforced by the competition rules. Although, such offline learning methods may be used alongside our approach to get even better results.
	
	This paper is organized as follows. Section \ref{sec:benchmarking} introduces the preliminary notations and the benchmarking system used in this paper. Section \ref{sec:primalsol} describes techniques to improve the primal bound by reusing the solutions obtained for the previous instances. Sections \ref{sec:branchinginfo} to \ref{sec:turnoff} present techniques to improve the dual bound. Section \ref{sec:branchinginfo} focuses on generating smaller branch-and-bound trees by reusing branching information. Section \ref{sec:paramtuning} describes a framework for automating the parameter tuning process. Section \ref{sec:turnoff} focuses on exploring the branch-and-bound trees faster by turning off inefficient solver parts. Section \ref{sec:finalcomparison} analyzes the performance of each technique and all of them together. Finally, Section \ref{sec:conclusion} summarizes our approach and presents future research directions for reoptimization. 
	
	\section{Benchmarking}
	\label{sec:benchmarking}
	\begin{comment}
		
		\begin{itemize}
			\item Hardware and software specification.
			\item 50 instances in 5 batches.
			\item Score computation. (time and gap only)
		\end{itemize}
		content...
	\end{comment}

	%\subsection{{Computational resources}}
	{\textbf{Computational Resources}}
	
	In this section, we discuss our computational experimental setup. We performed all the evaluations using SCIP optimization suite 8.0.2 \cite{BestuzhevaEtal2021ZR} that includes SOPLEX 6.0.2 on a CPU server with two sockets of Intel(R) Xeon(R) Gold 6258R CPU @ 2.70GHz 28 cores each (total of 56 cores) and 512 G{B} of RAM. We used only one thread for each instance and a maximum of 16 GB RAM. 
	
	%\subsection{{Dataset}}
	\noindent
	{\textbf{Dataset}}
	
	{The dataset described in this section will be used in all the experiments presented in this paper. These instances were used in the} MIP Computational Competition 2023. There are 15 series. Each series has 50 similar instances with one or more components changing across instances. {For each series, we solve the instances in the provided sequence, reusing only the information generated by solving previous instances in the series.} 
	
	\begin{comment}

	Consider a mixed integer program 
	
	\begin{align}
		\begin{split}
			{m}in &\quad c^Tx\\
			s.t. &\quad Ax \leq b\\
			&\quad l \leq x \leq u\\
			&\quad x_j \in \mathbb{Z}, \forall j \in I
		\end{split}
	\end{align}
	
	\noindent
	with $A\in \mathbb{R}^{m \times n}, c \in \mathbb{R}^n, b \in \mathbb{R}^m, l \in (\mathbb{R} \cup \{ -\infty \} )^n, u \in (\mathbb{R} \cup \{ \infty \} )^n$, variables $x \in \mathbb{R}^n$ with $I \subseteq N = \{1,...,n\}$.
	
			content...
\end{comment}
	
	The names of the series suggest the changing components. For each series, we used a time limit for solving each instance in that series. Table \ref{tab:seriescctl} describes the changing components{, time limit (in seconds), and shifted geometric means (shifted by 10 seconds) of solving times without using any reoptimization techniques for each series.\footnote{For the series RHS 2 and RHS 4, the shifted geometric means of solving times are higher than the time limit. This indicates that SCIP could not solve most of these instances in the given time limit.}}
	
	\begin{table}[h]
		\centering
		\caption{Changing components and time limits of different series.}
		\label{tab:seriescctl}
		\begin{tabular}{ll>{\raggedleft\arraybackslash}p{0.12\textwidth}r}
			\hline\noalign{\smallskip}
			Series & Changing component(s) & Time limit per instance & {Solving times}\\
			\noalign{\smallskip}\hline\noalign{\smallskip}
			Bound 1 & $l,u$ & 600 & 400.48 \\
			Bound 2 & $l,u$ & 300 & 186.23\\
			Bound 3 & $l,u$ & 600 & 336.96\\
			RHS 1 & $b$ & 400 & 294.39\\
			RHS 2 & $b$ & 60 & 69.55\\
			RHS 3 & $b$ & 550 & 307.66\\
			RHS 4 & $b$ & 60 & 70.00\\
			Objective 1 & $c$ & 400 & 301.95 \\
			Objective 2 & $c$ & 300 & 256.69\\
			Objective 3 & $c$ & 350 & 153.85\\
			RHS Objective 1 & $b,c$ & 500 & 492.80\\
			RHS Objective 2 & $b,c$ & 250 & 91.03\\
			Matrix 1 & $A$ & 300 & 270.70\\
			Matrix RHS Objective Bound & $A, b, c, l, u$ & 180 & 92.80\\
			Matrix RHS Bound & $A, b, l, u$ & 400 & 204.31\\
			\noalign{\smallskip}\hline
		\end{tabular}
	\end{table}

	{The detailed description of the origin of each series is provided in \cite{bolusani2023mip}. The series Bound 1, Bound 2, Bound 3, Objective 2, and RHS 3 are generated by mutating some instances from the MIPLIB 2017 benchmark library \cite{gleixner2021miplib}. Other series are based on various theoretical approaches (including stochastic multiple binary knapsack problem \cite{angulo2016improving}, column-generation-based heuristic \cite{firat2020column}, and constraint generation for MIP \cite{jimenez2022warm}) and real-world applications (including hydro unit commitment problem \cite{thomopulos2023generating}, multilevel supply chain \cite{sapse}, a hydroelectric valley industrial use case \cite{andreassian2021camels}, and optimal vaccine allocation problem \cite{tanner2010iis}).}

	\noindent
	{\textbf{Scoring}}
	
	For each instance in a series, we compute a `Total score' to measure the solver's performance on that instance, as used in the MIP Computational Competition 2023. The `Total score' {is the sum of} two components with equal weights: `Time score' and `Gap score'. 
	
	The `Time score' is the fraction $\dfrac{\text{time spent}}{\text{time limit}}$ if the instance is solved to optimality and $1$ otherwise. The `Gap score' is defined as the fraction $\dfrac{|pb - db|}{\text{max}(|pb|,|db|)}$, where $pb$ and $db$ are the primal and dual bounds, respectively. The `Gap score' is $1$ if the primal or dual bound are infinite or have different signs, i.e., $pb \cdot db < 0$. {When the instance is solved to optimality, the `Gap score' is 0, and when the `Gap score' is positive, the `Time score' is 1 as the solver failed to prove optimality in the given time limit.}
	
	The competition also included a `Nofeas score' in the total score, which is equal to $1$ if the solver failed to find any feasible solution for the instance {(implying `Time score' and `Gap score' to be 1)}. In all our experiments, the solver always found at least one feasible solution {for each instance}, so this component was $0$ for all instances; hence, we removed it {for simplicity}.
	
	For each experiment, we report the average instance total score for batches of 10 instances (1-10, 11-20, 21-30, 31-40, 41-50). The average total score for each batch helps us understand the effect of the reoptimization techniques as we progress through the series.
	
	\section{Reusing primal information}
	\label{sec:primalsol}
	\begin{comment}
		
		\begin{itemize}
			\item {Increase the efforts of \textit{completesol} heuristic.}
			\item {Using previous solution as it is.}
			\item {Removing the continuous variables + clipping.}
			\item {Common solution}
			\item {Conversion numbers}
			\item Comparison with base.
		\end{itemize}
		content...
	\end{comment}
	
	When solving multiple similar instances, we expect the solutions of those instances to be similar as well. In this section, we describe various ways to use the solutions of the instances solved before.
	
	{Modern MIP solvers aim to derive a feasible solution from any complete or partial solution provided to them (denoted as a `solution hint' or `hint' henceforth) at the beginning of the search process. In SCIP, this is done by a heuristic called \textit{completesol}. SCIP can take multiple solution hints at the beginning of the search process. The \textit{completesol} heuristic processes each hint sequentially. The heuristic checks for feasibility if the hint represents a complete solution. If feasible, the solution is added to the solution pool when it improves the incumbent. For a partial solution, the heuristic attempts to complete it by solving a sub-MIP. To our knowledge, the heuristic does not attempt to repair an infeasible solution hint. 
		
	Unless only the objective function has changed (i.e., all feasible solutions of the previous instances are also feasible for the current instance), we update the node limit (number of nodes it can use in the sub-MIP) for this heuristic to 5000 (default limit = 500) and also remove the limit on the maximum number of improving solutions (default limit = 5). These changes ensure that the \textit{completesol} spends more effort to convert a given solution hint into a feasible solution.}
	
	{Removing the limit on the maximum number of improving solutions allows us to provide as many solution hints as possible. However, we expect the best-found solutions of the previous 10 instances to represent well the solution pattern (or trend) of the series. Hence, we decided not to use more than 10 solution hints in any series (5 for the Objective series) for our experiments. If the solution hints help improve the total score, one can experiment with providing more solution hints.}
	
	First, we experimented with providing the complete {best-found} solutions for the {previous} 10 instances (if available). For series where only the objective has changed, we provide the complete solutions for the {previous} 5 instances to reduce the efforts spent by the \textit{completesol} heuristic.
	
	{Next, w}e remove the continuous variables' values from the hint {(making the solution a partial solution)} to improve the {number of hints getting converted into a feasible solution (hint-to-solution conversion number)}. 
	
	Furthermore, {when the variable bounds change in a series,} we clip the solution values by the bounds of the current instance. If a variable's value in the hint is greater than its new upper bound, we set the value to the new upper bound in the hint. We use similar clipping for the lower bound. 
	
	{For this experiment, a}s before, we provide clipped integer solution hints for the {previous} 10 instances except in the Objective series, where we provide solution hints for the {previous} 5 instances. If the assignment of all integer variables can produce a feasible solution, the \textit{completesol} heuristic can find it by solving a linear program.
	
	{Finally, w}e construct a common partial solution to improve the hint-to-solution conversion number {as follows}. We record the best{-found} solutions of the instances. We add a variable-value pair to the hint that is part of the first instance solution and present in at least $\alpha$\% of the {(best-found)} solutions {of all previous instances}. {Being part of the first instance solution} ensures that the provided hint is feasible for at least one instance {in the series} and reduces the chances of combining inconsistent variable-value pairs for the entire series. This idea {of combining multiple solutions} is inspired by many existing heuristics like RINS \cite{danna2005rins}, Crossover \cite{berthold2006crossover}, and CMSA \cite{blum2016cmsa}. 
	
	{The $\alpha$ value controls the amount of effort needed to find a feasible solution (using the given solution hint) and the solution's quality. After preliminary testing of different values of $\alpha \in \{80, 90, 95, 100\}$, we selected the value $\alpha = 90$ for our experiments.} We did not experiment with dynamically changing $\alpha$ for a series. 
	
	{For our final experiment, w}e provide a common partial solution hint and the clipped integer solutions of the {previous} 9 instances (4 instances for the Objective series). 
	
	To measure the effect of each experiment, we record the percentage of instances for which the provided hints are converted into at least one feasible solution by the \textit{completesol} heuristic. The conversion percentages are reported in Table \ref{tab:hintconv}.
	
	{C}omplete solution hints are converted into feasible solutions for all instances in all Objective series. The conversion percentages are above 80\% for both RHS Objective series. For the other series, the conversion numbers are very low, as expected. 
	
	Removing continuous variable values from the hint {along with bound clipping} slightly improves the conversion percentages. {The bound clipping only affects five series where the variable bounds are changing. The improvement in the conversion percentage is largely due to the removal of continuous variables' values.} 
	
	Providing a common solution hint significantly improves the conversion percentages. However, the conversion percentages remain 0 for RHS series 3 and Matrix RHS Bound series. 
	
	\begin{table}[h]
		\centering
		\caption{Percentages of instances for which the provided hints are converted into at least one feasible solution by \textit{completesol} heuristic for different types of hints. {Better numbers are in bold.}}
		\label{tab:hintconv}
		\begin{tabular}{p{0.4\textwidth}
				>{\raggedleft\arraybackslash}p{0.12\textwidth}
				>{\raggedleft\arraybackslash}p{0.12\textwidth}
				>{\raggedleft\arraybackslash}p{0.12\textwidth}}
			\hline\noalign{\smallskip}
			Series & Complete solution (\%) & Clipped Integer solution (\%) & Common Clipped Integer solution (\%) \\
			\noalign{\smallskip}\hline\noalign{\smallskip}
			Bound 1 & 40.8 & 40.8 & \textbf{51.0} \\
			Bound 2 & 4.1 & 6.1 & \textbf{42.9} \\
			Bound 3 & 8.2 & 10.2 & \textbf{61.2} \\
			RHS 1 & 0.0 & 0.0 & \textbf{69.4} \\
			RHS 2 & \textbf{100.0} & \textbf{100.0} & \textbf{100.0} \\
			RHS 3 & 0.0 & 0.0 & 0.0 \\
			RHS 4 & \textbf{100.0} & \textbf{100.0} & \textbf{100.0} \\
			Objective 1 & \textbf{100.0} & \textbf{100.0} & \textbf{100.0} \\
			Objective 2 & \textbf{100.0} & \textbf{100.0} & \textbf{100.0} \\
			Objective 3 & \textbf{100.0} & \textbf{100.0} & \textbf{100.0} \\
			RHS Objective 1 & 81.6 & 85.7 & \textbf{93.9} \\
			RHS Objective 2 & 83.7 & \textbf{89.8} & \textbf{89.8} \\
			Matrix 1 & 0.0 & \textbf{100.0} & \textbf{100.0} \\
			Matrix RHS Objective Bound & 65.3 & 65.3 & \textbf{87.8}\\
			Matrix RHS Bound &0.0 & 0.0 & 0.0 \\
			\noalign{\smallskip}\hline\noalign{\smallskip}
			Average &52.2	&59.9 &	\textbf{73.1}\\
			\noalign{\smallskip}\hline
		\end{tabular}
	\end{table}

	Finally, we show the effect of providing solution hints on the total scores for each series in Table \ref{tab:solhintaggr}. {The conversion rate column shows the percentage of instances for which the solution hints were converted into at least one feasible solution by the solver (same as the last column in Table \ref{tab:hintconv}).} The BASE column shows the total scores of the variant where we do not provide any solution hints{, and the SOLHINT column shows the percentage improvement over BASE scores for} the variant where we provide the Common Clipped Integer solution hints as described above. {The detailed batch-wise results are presented in Table \ref{tab:solutionhintscores}.}
	
	In the series Bound 3, {RHS 1}, RHS 2, RHS 4, Objective 1, {RHS Objective 1}, and Matrix RHS Objective Bound, the conversion {rates} are above 50\%, and we observe a gain in the total scores. For the other series, it is unclear whether providing solution hints is helpful. Specifically for series Bound 1 and Matrix 1, providing solution hints significantly degrades the total scores even though the conversion rate is high. We address this issue in Section \ref{sec:paramtuning}.

\begin{table}[h]
	\centering
	\caption{{Effects of providing solution hints on total scores. Average BASE scores and \% improvements caused by SOLHINT variant. Positive improvement numbers are in bold.}}
	\label{tab:solhintaggr}
	\begin{tabular}{p{0.4\textwidth}
			>{\raggedleft\arraybackslash}p{0.12\textwidth}
			>{\raggedleft\arraybackslash}p{0.12\textwidth}
			>{\raggedleft\arraybackslash}p{0.15\textwidth}}
		Series & Conversion rate (\%) & BASE & SOLHINT (\% improvement over BASE)\\
		\hline\noalign{\smallskip}
		Bound 1 & 51.0 & 0.734 & -10.90 \\
		Bound 2 & 42.9 & 0.658 & -0.30 \\
		Bound 3 & 61.2 & 0.636 & \textbf{2.52} \\
		RHS 1 & 69.4 & 0.746 & \textbf{2.41} \\
		RHS 2 & 100.0 & 0.994 & \textbf{2.82} \\
		RHS 3 & 0.0 & 0.652 & -0.31 \\
		RHS 4 & 100.0 & 1.000 & \textbf{1.60} \\
		Objective 1 & 100.0 & 0.756 & \textbf{4.50} \\
		Objective 2 & 100.0 & 0.910 & -3.52 \\
		Objective 3 & 100.0 & 0.420 & -1.43 \\
		RHS Objective 1 & 93.9 & 0.996 & \textbf{4.42} \\
		RHS Objective 2 &89.8 & 0.418 & \textbf{0.48} \\
		Matrix 1 & 100.0 & 0.896 & -6.47 \\
		Matrix RHS Objective Bound &87.8 & 0.504 & \textbf{28.17} \\
		Matrix RHS Bound & 0.0 & 0.572 & 0.00\\
		\noalign{\smallskip}\hline\noalign{\smallskip}
		Average & 73.1 & 0.726 & \textbf{1.60} \\
		\noalign{\smallskip}\hline
	\end{tabular}
\end{table}
	
	\section{Reusing and improving branching information}
	\label{sec:branchinginfo}
	
	\subsection{Reusing branching information}
	\label{sec:reusepscost}
	\begin{comment}
		
		\begin{itemize}
			\item {History and global history}
			\item{Reduced count of pseudocost.}
			\item Comparison with base.
		\end{itemize}
		
	\end{comment}
	
	MIP solvers spend a significant amount of time proving the solution's optimality. The MIP solvers record various branching-related statistics in the branch and bound process{ that are used} in branching heuristics.
	
	SCIP stores the branching statistics for each variable in the variable history data structure. This structure stores many branching-related records, including the variable's pseudocost information, conflict records, and inference records. 
	
	For each instance in a series, we retrieve each variable's history. We pass this variable history to the corresponding variable in the next instance at the beginning of the solving process. Since the instances are similar, we expect this information to also be {useful} for the next instance{, and t}he solver no longer needs to compute the branching statistics from scratch. 
	
	Along with branching statistics for each variable, SCIP maintains a `global history' structure containing the {same} statistics for the entire problem {(i.e., statistics for all variables combined)}. {When the number of branching records for a particular variable is very low, SCIP uses the records stored in global history for that variable while making branching decisions.} We also extract this {global history} structure and pass it to the next instance.
	
	{The current state-of-the-art branching rules heavily rely on the pseudocosts for variables \cite{achterberg2005reliability}. Pseudocosts measure the average objective improvement for a unit change in value for every integer variable. Their use in branching was first reported in \cite{benichou1971pscost}.} For each variable, the solver maintains the pseudocost and updates it incrementally by storing the number of updates as {the} pseudocost count in the history structure.
	
	One of the most commonly used branching rules is reliable pseudocosts {\cite{achterberg2005reliability}}. At the beginning of the solving process, the pseudocosts {(and other records in the variable history)} are unreliable {since the number of pseudocost updates is low \cite{linderoth1999searchstrategies}}. To address this, the reliable pseudocosts branching rule uses strong branching {\cite{applegate1995strongbranching}} at the beginning of the solving process. In strong branching, the solver evaluates multiple variables for branching at each node and updates their pseudocosts. Hence, even though the branching rule takes longer to execute, the variables receive more pseudocost updates. This is required to have reliable pseudocosts.
	
	 When we pass the history to the next instance, we reduce the pseudocost count to 4 if it is greater than 4. The idea behind this reduction is to make the solver use the pseudocosts of the previous instance as a `warm start' only. This allows the pseudocosts to converge to the actual average value (which can be different for the next instance) faster. This is inspired by \cite{berthold2022restarttransfer}. We used the value 4 for reducing the count as a threshold because it is slightly less than the fixed number threshold 5, after which pseudocosts are considered `reliable' in SCIP \cite{hendel2015pscostvariance}. {We expect such a warm start to only reduce the computing time by reducing the number of strong branching calls and not the number of branching nodes explored.}
	
	\subsection{Improving pseudocosts}
	\label{sec:strongbranching}
	\begin{comment}
		\begin{itemize}
			\item Idea from reliability branching.
			\item Comparison with base and base + reuse branching.
		\end{itemize}
	\end{comment}

	We extend the reliability branching idea further and solve the first instance in the series with full strong branching to generate better pseudocosts. Since we reuse the pseudocosts across different instances, generating more reliable pseudocosts is essential. Solving the instance with strong branching helps us generate more reliable pseudocosts.
	
	Furthermore, we observed that the pseudocosts are more reliable when the objective and variable bounds are unchanged. In that case, we switch to the pseudocost branching rule instead of the default reliable pseudocost branching rule after solving the first instance with strong branching.
	
	Using full strong branching comes at a cost. The solving process is significantly slower for the first instance, where we use full strong branching. {In only three series (compared to nine in BASE approach), we could solve the first instance to optimality using strong branching: RHS Objective 2, Matrix 1, and Matrix RHS Bound. The average gap score for the first instance across all series is 0.18 when we use strong branching and 0.003 when we do not.} 
	
	However, we expect the recorded pseudocosts to be more reliable and compensate for the loss on the first instance by solving subsequent instances faster. {In the competition, the total score for the first instance in a given series was given smaller weights compared to the subsequent instances, motivating us to experiment with more exploration for the initial few instances.}
	
	Table \ref{tab:reusehistaggr} shows the effect of reusing branching information and using strong branching for the first instance on the total scores. The BASE column shows the total scores of the variant where we do not reuse any branching information{, t}he REHIS column shows the {percentage improvement over BASE }total scores of the variant where we reuse the variable and global history{, and} the column REHIS+SB shows the {percentage improvement over BASE }total scores for the variant where we reuse the variable and global history and solve the first instance with strong branching. {The detailed batch-wise results are presented in Table \ref{tab:reusehistoryscores}.}
	
	The scores in Table \ref{tab:reusehistaggr} suggest that reusing branching history improves the solver's performance in general. However, it is not better when the objective is changing. This is expected because we should not expect the pseudocosts to remain similar across instances where the objective changes. Using strong branching for the first instance is better on all RHS series. {We do not have a good explanation for why using strong branching is most effective on series with RHS changing.} {F}or {the Bound and the Matrix} series, it does not show a clear improvement in the scores over the REHIS variant. 
	
	\begin{table}[h]
		\centering
		\caption{{Effects of reusing branching history on total scores. Average BASE scores and \% improvements caused by REHIS and REHIS+SB variants. Better positive improvement numbers are in bold.}}
		\label{tab:reusehistaggr}
		\begin{tabular}{p{0.4\textwidth}
				>{\raggedleft\arraybackslash}p{0.12\textwidth}
				>{\raggedleft\arraybackslash}p{0.15\textwidth}
				>{\raggedleft\arraybackslash}p{0.15\textwidth}}
			Series & BASE & REHIS (\% improvement over BASE) & REHIS+SB(\% improvement over BASE) \\
			\hline\noalign{\smallskip}
			Bound 1 & 0.734 & \textbf{34.06} & 30.25 \\
			Bound 2 & 0.658 & \textbf{12.77} & 6.08 \\
			Bound 3 & 0.636 & \textbf{4.72} & 3.14 \\
			RHS 1 & 0.746 & 6.43 & \textbf{7.77} \\
			RHS 2 & 0.994 & 0.00 & \textbf{6.64} \\
			RHS 3 & 0.652 & 7.98 & \textbf{14.42} \\
			RHS 4 & 1 & 0.00 & \textbf{5.80} \\
			Objective 1 & 0.756 & 2.38 & \textbf{3.17} \\
			Objective 2 & 0.91 & -6.15 & -11.21 \\
			Objective 3 & 0.42 & -12.38 & -15.71 \\
			RHS Objective 1 & 0.996 & \textbf{0.20} & -1.61 \\
			RHS Objective 2 & 0.418 & -9.57 & -3.35 \\
			Matrix 1 & 0.896 & -1.79 & -4.25 \\
			Matrix RHS Objective Bound & 0.504 & -19.44 & -17.06 \\
			Matrix RHS Bound & 0.572 & -0.69 & \textbf{1.39} \\
			\noalign{\smallskip}\hline\noalign{\smallskip}
			Average & 0.726 & 1.23 & \textbf{1.70} \\
			\noalign{\smallskip}\hline
		\end{tabular}
	\end{table}

	\section{Automated parameter tuning}
	\label{sec:paramtuning}
	\begin{comment}
		
		\begin{itemize}
			\item UCB algorithm equation modifications.
			\item Scores and C value.
			\item Using standard deviation.
			\item Forced initial exploration.
			\item Tuning sequence.
			\item Effect of Cuts tuning. (both together). And stats about chosen params.
			\item Effect of primal hint tuning. And stats about chosen params.
		\end{itemize}
	\end{comment}
	
	MIP solvers contain many tunable parameters{, and s}ince we are solving similar instances in a series, we want to automate the parameter tuning process to improve the solving speed for the later instances. {The parameter tuning with offline training has been explored before (see, e.g., \cite{xu2011hydra}). This section focuses on MIP solver parameter tuning without any offline training.}
	
	For a given parameter, {we consider two values for tuning (i.e., ON and OFF) and maintain scores for those values.} The score computation is inspired by the Upper Confidence Bound (UCB) Algorithm {\cite{auer2002finite}} from reinforcement learning. We compute the score $S_v$ for a value $v$ using equation (\ref{eq:param_score}). 
		
	\begin{equation}
		\label{eq:param_score}
		S_v = Q_v + C / N_v
	\end{equation}
	
	It involves two parts: ({a}) The running average of the base score ($Q_v$), and ({b}) The uncertainty of the score, i.e., the confidence bound ($1/N_v$), where $N_v$ is the number of score updates the value has received. When a value is not explored enough, the confidence bound adds a higher number to the score to encourage more exploration for that value. $C$ is a weight constant that determines how fast the score converges. Based on our experiments, the value $C=0.3$ gives the best results regarding convergence correctness {(i.e., the desired value being used more)} and {convergence} speed. For the base score, we use the negative total score of the instance, as described in Section \ref{sec:benchmarking}.
	
	In the UCB algorithm, the confidence bound uses the square root of $N_v$ and the total number of score updates across all values. Using the total number of score updates ensures that confidence in the score of a value decreases over time. These parts are helpful in general. However, we discard the total number of score updates for the competition and use $N_v$ instead of its square root for faster convergence. {See Table \ref{tab:convergence} for a comparison of convergence speed. Using $\sqrt{N_v}$ results in more exploration even in the last batch (instances 41-50) for multiple values of $C$, which is not desired for the competition but may be useful in practice when the total number of instances is significantly more than 50.}
	
	\begin{table}[H]
		\centering
		\caption{{Convergence comparison between $C/N_v$ and $C/\sqrt{N_v}$.}}
		\label{tab:convergence}
		\begin{tabular}{rrrrrrr}
			\hline\noalign{\smallskip}
			& \multicolumn{2}{c}{C=0.2} & \multicolumn{2}{c}{C=0.3} & \multicolumn{2}{c}{C=0.4} \\
			\noalign{\smallskip}\cline{2-7}\noalign{\smallskip}
			$N_v$ & $C/N_v$ & $C/\sqrt{N_v}$ & $C/N_v$ & $C/\sqrt{N_v}$ & $C/N_v$ & $C/\sqrt{N_v}$ \\
			\noalign{\smallskip}\hline\noalign{\smallskip}
			1 & 0.200 & 0.200 & 0.300 & 0.300 & 0.400 & 0.400 \\
			%2 & 0.100 & 0.141 & 0.150 & 0.212 & 0.200 & 0.283 \\
			%3 & 0.067 & 0.115 & 0.100 & 0.173 & 0.133 & 0.231 \\
			%4 & 0.050 & 0.100 & 0.075 & 0.150 & 0.100 & 0.200 \\
			5 & 0.040 & 0.089 & 0.060 & 0.134 & 0.080 & 0.179 \\
			10 & 0.020 & 0.063 & 0.030 & 0.095 & 0.040 & 0.126 \\
			%15 & 0.013 & 0.052 & 0.020 & 0.077 & 0.027 & 0.103 \\
			20 & 0.010 & 0.045 & 0.015 & 0.067 & 0.020 & 0.089 \\
			%25 & 0.008 & 0.040 & 0.012 & 0.060 & 0.016 & 0.080 \\
			30 & 0.007 & 0.037 & 0.010 & 0.055 & 0.013 & 0.073 \\
			%35 & 0.006 & 0.034 & 0.009 & 0.051 & 0.011 & 0.068 \\
			40 & 0.005 & 0.032 & 0.008 & 0.047 & 0.010 & 0.063 \\
			%45 & 0.004 & 0.030 & 0.007 & 0.045 & 0.009 & 0.060 \\
			50 & 0.004 & 0.028 & 0.006 & 0.042 & 0.008 & 0.057 \\
			100 & 0.002 & 0.020 & 0.003 & 0.030 & 0.004 & 0.040 \\
			200 & 0.001 & 0.014 & 0.002 & 0.021 & 0.002 & 0.028 \\
			500 & 0.000 & 0.009 & 0.001 & 0.013 & 0.001 & 0.018 \\
			\noalign{\smallskip}\hline
		\end{tabular}
	\end{table}

	The tuning involves exploration and exploitation stages. The score value $S_v$ already considers the number of score updates and motivates some exploration. In addition, we mark a value under exploration if it has not been used at least four times. For a given parameter, we prioritize selecting the values under exploration, if any.
	
	For the exploitation stage, selecting the value with the highest score $S_v$ is sufficient. However, the UCB algorithm assumes that the underlying reward distribution is constant. This would be the case if all instances have similar difficulties. However, this is not always true. The instance difficulties may change in a series. To address that, we also use the standard deviation of the base score. For the selection, we mark any values $v$ with a score $S_v$ within a distance $1/10$ of the standard deviation from the highest score $S_v^*$ as candidate values. If there are multiple candidate values, we select the final value {uniformly at random} among them.
	
	In this framework, one can tune multiple parameters at the same time. However, doing so increases the noise, and the scores converge slowly. We tune three binary parameters in our approach.
	
	\begin{enumerate}
		\item Provide hint or not
		\item Use root node cuts or not
		\item Use cuts at other nodes or not
	\end{enumerate}
	
	In general, providing hints helps the solver. However, as shown in Section \ref{sec:primalsol}, in some cases, providing solution hints causes harm in the solving process. Hence, we decided to use our parameter tuning framework to decide whether to provide solution hints to the solver. For this parameter, we update the score for value `ON' (provide hint) only if the solver converts the hints into at least one feasible solution. Otherwise, we update the value `OFF' score if we do not provide hints or the solver fails to convert the hints into a feasible primal solution.
	
	{Similarly, in most cases, using cuts is helpful. However, adding cuts can hurt the performance in some instances. We expect the effects of adding cuts to remain similar throughout the entire series, i.e., it either improves or degrades performance for all instances. But, we expect the effect to differ between root node cuts and other node cuts. The root node cuts are always global (i.e., valid for the entire tree), and solvers typically spend more effort generating cuts at the root node. For the cuts parameters above, we do not check if the solver adds the cuts when the `ON' value is selected, unlike the solution hint parameter.}
	
	When tuning multiple parameters, we perform initial exploration in a deterministic fashion. The first parameter selects the value `ON' on every alternate instance during exploration. The second parameter selects the value `ON' on instances with the second bit set on the binary representation of the instance index. The third parameter selects the value `ON' on instances with the third bit set on the binary representation of the instance index. {Table \ref{tab:detexploration} shows how the values are selected for each parameter during deterministic exploration if tuning of all parameters starts at instance 0.}
	
		\begin{table}[H]
		\centering
		\caption{{Selected parameter values during deterministic exploration.}}
		\label{tab:detexploration}
		\begin{tabular}{clll}
			\hline\noalign{\smallskip}
			Instance index & Solution hint & Cuts & Root cuts  \\
			\noalign{\smallskip}\hline\noalign{\smallskip}
			0 & OFF & OFF & OFF \\
			1 & ON & OFF & OFF \\
			2 & OFF & ON & OFF \\
			3 & ON & ON & OFF \\
			4 & OFF & OFF & ON \\
			5 & ON & OFF & ON \\
			6 & OFF & ON & ON \\
			7 & ON & ON & ON \\
			\noalign{\smallskip}\hline
		\end{tabular}
	\end{table}
	
	For the computational experiments, we measure the effects of parameter tuning in the presence of other reoptimization techniques discussed in Sections \ref{sec:primalsol} {and} \ref{sec:branchinginfo}. The parameters' tuned values may differ if the other reoptimization techniques are not used. We noticed that disabling restarts gives better performance in most cases. The probable cause might be that the solver is already reusing much information from the previous solutions. Hence, for the computational experiments, we turn off restarts. 
	
	Table \ref{tab:paramstats} shows the most selected value for each parameter for each series. The count column shows the number of times that value was used out of 49. We do not count the parameter values used for the first instance. The count values suggest how quickly the above system converged for the selected values. Some parameters do not significantly affect the total scores and did not converge. For example, the selected values for the solution hint parameter in Bound 2, Objective 2, and Objective 3 series have counts close to 25. {The count values above 40 in most cases suggest a strong preference, and count values between 30 and 40 suggest a weak preference for a parameter value. Note that these preferences may not always be accurate (i.e., we get better scores when we use the preferred values.) }
	
	{Because of the initial deterministic exploration, we expect the max count value to be 45. However, the count values for the solution hint parameter are 49 for the RHS 3 and Matrix RHS Bound series. This happens because the count for the `ON' value in this parameter does not represent how many times the hints were provided but how many times the hints were converted into a feasible solution. None of the solution hints were converted into a feasible solution for the RHS 3 and Matrix RHS Bound series (see Table \ref{tab:hintconv}). Hence, even when the hints were provided, the scores of the `OFF' values were updated. Since the counts for the `ON' value never reached count 4, it remains marked under exploration, and the tuning algorithm keeps using it.}
	
	\begin{table}[H]
		\centering
		\caption{Selected values and count for the tuned parameters for each series.}
		\label{tab:paramstats}
		\begin{tabular}{llrlrlr}
			\hline\noalign{\smallskip}
			Series & \multicolumn{2}{c}{Solution hint} & \multicolumn{2}{c}{Cuts} & \multicolumn{2}{c}{Root cuts} \\
			\noalign{\smallskip}\cline{2-7}\noalign{\smallskip}
			& Value & Count & Value & Count & Value & Count \\
			\noalign{\smallskip}\hline\noalign{\smallskip}
			Bound 1 & OFF & 45 & ON & 45 & ON & 45 \\
			Bound 2 & OFF & 27 & ON & 45 & ON & 45 \\
			Bound 3 & OFF & 45 & ON & 45 & ON & 45 \\
			RHS 1 & ON & 33 & ON & 45 & OFF & 41 \\
			RHS 2 & ON & 45 & OFF & 41 & OFF & 41 \\
			RHS 3 & OFF & 49 & ON & 40 & OFF & 34 \\
			RHS 4 & ON & 45 & OFF & 41 & OFF & 41 \\
			Objective 1 & OFF & 28 & OFF & 32 & OFF & 29 \\
			Objective 2 & ON & 27 & ON & 43 & ON & 37 \\
			Objective 3 & ON & 39 & ON & 45 & OFF & 41 \\
			RHS Objective 1 & ON & 30 & OFF & 41 & OFF & 41 \\
			RHS Objective 2 & ON & 36 & OFF & 28 & ON & 45 \\
			Matrix 1 & ON & 38 & ON & 45 & ON & 44 \\
			Matrix RHS Objective Bound & ON & 40 & ON & 45 & OFF & 41 \\
			Matrix RHS Bound & OFF & 49 & ON & 45 & ON & 45\\
			\noalign{\smallskip}\hline
		\end{tabular}
	\end{table}
	
	Table \ref{tab:paramtuningaggr} shows the effect of automated parameter tuning on the total scores. The column NOTUNING shows the total scores of the variant where we do not tune the parameters{, and t}he column TUNING shows the total scores of the variant where we automate the parameter tuning as described above. For both variants, we provide the solution hints as described in Section \ref{sec:primalsol} and reuse the branching information as described in Section \ref{sec:branchinginfo}. For both variants, we turn off restarts in the solver. {The detailed batch-wise results are presented in Table \ref{tab:paramtuningscores}.}
	
	In most series, the variant TUNING performs worse in the first batch. This is because of the exploration stage of the tuning process. The TUNING variant scores for the later batches are better compared to the NOTUNING variant in general. For the series Bound 3, RHS Objective 2, and Matrix 1, the TUNING variant has worse scores than the NOTUNING variant {because} the TUNING variant converged to incorrect values for some parameters {(see the batch-wise comparison in Table \ref{tab:paramtuningscores})}. In the Objective 2{ and the RHS 3} series, it is unclear if the TUNING variant converged to the correct values. {In the series Bound 2 and Matrix RHS Bound, the TUNING variant improves the total score in the subsequent batches, but it does not cover the losses made in the exploration stage. In the other series, the TUNING variant clearly outperforms the NOTUNING variant.}

		\begin{table}[h]
		\centering
		\caption{{Effects of automated parameter tuning on total scores. Average NOTUNING scores and \% improvements caused by TUNING variant. Positive improvement numbers are in bold.}}
		\label{tab:paramtuningaggr}
		\begin{tabular}{p{0.4\textwidth}
				>{\raggedleft\arraybackslash}p{0.12\textwidth}
				>{\raggedleft\arraybackslash}p{0.15\textwidth}}
			Series & NOTUNING & TUNING (\% improvement over NOTUNING)  \\
			\hline\noalign{\smallskip}
			Bound 1 & 0.7 & \textbf{26.86} \\
			Bound 2 & 0.444 & -3.60 \\
			Bound 3 & 0.484 & -14.05 \\
			RHS 1 & 0.672 & \textbf{3.87} \\
			RHS 2 & 0.664 & \textbf{12.35} \\
			RHS 3 & 0.596 & -5.03 \\
			RHS 4 & 0.678 & \textbf{13.57} \\
			Objective 1 & 0.744 & \textbf{11.02} \\
			Objective 2 & 1.008 & -2.18 \\
			Objective 3 & 0.656 & \textbf{38.72} \\
			RHS Objective 1 & 0.992 & \textbf{6.05} \\
			RHS Objective 2 & 0.572 & -27.62 \\
			Matrix 1 & 1.016 & -2.56 \\
			Matrix RHS Objective Bound & 0.466 & \textbf{14.16} \\
			Matrix RHS Bound & 0.726 & -6.61 \\
			\noalign{\smallskip}\hline\noalign{\smallskip}
			Average & 0.695 & \textbf{4.33 }\\
			\noalign{\smallskip}\hline
		\end{tabular}
	\end{table}

	\section{Other optimizations}
	\label{sec:turnoff}
	\begin{comment}
		\begin{itemize}
			\item Unsuccessful Separators
			\item Unsuccessful Heuristics
			\item Unsuccessful presolve rules.
		\end{itemize}
	\end{comment}
	
	MIP solvers have many components turned on by default. Not all of them perform well during the solving process. After solving each instance, we extract the statistics for primal heuristics, presolvers, and separators (cut generators). If they are unsuccessful for the initial few instances and expensive in solving time, we turn them off to give the other components some extra time. Turning such non-performing components off had a minor improvement in the total scores. 
	
	We use the following criteria to mark a component unsuccessful. 
	
	\begin{itemize}
		\item A presolve rule is unsuccessful if the total number of changes made in the first 15 instances by that rule is 0.
		\item A separator is unsuccessful if the total number of generated cuts in the first 25 instances by that separator is 0.
		\item A primal heuristic is unsuccessful if it either fails to find any solution in the first 25 instances or the time spent per best solution found is more than 20\% of the given time limit for a single instance (measured after the first 25 instances).
	\end{itemize}
	
	We show the effect of turning inefficient solver parts off on total scores in Table \ref{tab:turnoffaggr}. {For this experiment, the other reoptimization techniques discussed in Sections \ref{sec:primalsol}, \ref{sec:branchinginfo}, and \ref{sec:paramtuning} were also in use,  ensuring that we only turn off a solver part if it is ineffective under various settings (different parameter values and improved pseudo costs). The detailed batch-wise results are presented in Table \ref{tab:turnoffscores}.} 
	
	Unlike previous experiments, turning the ineffective components does not consistently improve the total scores for most series. However, for series RHS 2, RHS 4, Objective 1, and RHS Objective 2, we see a significant improvement in the total scores of later batches. On the other hand, the total scores for RHS 3, {Objective 3, } and Bound 3 show significant degradation {in most batches}, suggesting that we turned off one or more essential components of the solver that were ineffective in the earlier instances of that series. {For the other series, it is unclear whether the detection of ineffective parts was successful (see the batch-wise comparison in Table \ref{tab:turnoffscores}).}
	
		\begin{table}[h]
		\centering
		\caption{{Effects of automated parameter tuning on total scores. Average NOTURNOFF scores and \% improvements caused by TURNOFF variant. Positive improvement numbers are in bold.}}
		\label{tab:turnoffaggr}
		\begin{tabular}{p{0.4\textwidth}
				>{\raggedleft\arraybackslash}p{0.15\textwidth}
				>{\raggedleft\arraybackslash}p{0.15\textwidth}}
			Series & NOTURNOFF & TURNOFF (\% improvement over NOTURNOFF)  \\
			\hline\noalign{\smallskip}
			Bound 1 & 0.508 & \textbf{5.51} \\
			Bound 2 & 0.434 & -5.07 \\
			Bound 3 & 0.52 & -2.69 \\
			RHS 1 & 0.66 & \textbf{3.64} \\
			RHS 2 & 0.57 & \textbf{24.56} \\
			RHS 3 & 0.584 & -16.10 \\
			RHS 4 & 0.582 & \textbf{21.65} \\
			Objective 1 & 0.648 & \textbf{2.78} \\
			Objective 2 & 1.046 & \textbf{1.72} \\
			Objective 3 & 0.398 & -8.04 \\
			RHS Objective 1 & 0.898 & -3.79 \\
			RHS Objective 2 & 0.686 & \textbf{5.54} \\
			Matrix 1 & 1.034 & -1.35 \\
			Matrix RHS Objective Bound & 0.36 & -2.22 \\
			Matrix RHS Bound & 0.766 & \textbf{0.52}\\
			\noalign{\smallskip}\hline\noalign{\smallskip}
			Average & 0.646 & \textbf{1.78 }\\
			\noalign{\smallskip}\hline
		\end{tabular}
	\end{table}

	\section{All together}
	\label{sec:finalcomparison}
	\begin{comment}
		\begin{itemize}
			\item Public instances
			\item Private instances
			\item Matrix changes
		\end{itemize}
	\end{comment}
	
	Table \ref{tab:finalscoresav} shows the average total scores of different techniques we used for reoptimization. The first row shows the average total score of solving each instance from scratch (Base). Reusing the primal solutions of previous instances has the least effect (1.30\%) of all the techniques described in this paper. Reusing branching history along with using strong branching for the first instance improves the performance slightly more (2.57\%). The most considerable improvement (of 4.53\% ) comes from the parameter tuning, which is more than using the primal solution hints and reusing branching history combined (4.40\%). Finally, combining all techniques results in a 12.77\% improvement over the Base approach.
	
	Furthermore, we show the performance of the reoptimization techniques per batch in Table \ref{tab:finalscoresbatchav}. Compared to the Base approach, our approach gives worse scores in the first batch by almost 4\%. {This slowdown is caused mainly by the parameter tuning exploration stage (see Table \ref{tab:paramtuningscores}). Using strong branching for the first instance also contributes to this slowdown but not as much as parameter tuning (see Table \ref{tab:reusehistoryscores}).}  However, our approach progressively strengthens and tunes the solver in the later batches. In the last three batches, the scores improve by more than 15\% compared to the Base approach. In the last batch, we observe an improvement of 22\% in the total scores compared to the Base approach.
	
	{We expect this improvement trend to slow down if we have more instances. We designed the approach to complete most solver tuning and improvements by the beginning of the last batch (41-50) and exploit what works well after that. However, this approach can be modified for real-world settings with more than 50 instances in a series. We talk more about such modifications in Section \ref{sec:conclusion}.}
	
	{Finally, we show the shifted geometric means of solving time (shifted by 10 seconds) for instances in each series in Table \ref{tab:stgeometricmeans}. Compared to the Base approach, our approach significantly improves the solving times for all series except RHS Objective 2, Matrix 1, and Matrix RHS Bound. In RHS Objective 2 and Matrix RHS Bound series, the biggest degradation is caused by automated parameter tuning. In Matrix 1, all of the techniques presented above degraded the performance. The biggest degradation is caused by reusing the primal solutions, which the automated parameter tuning failed to stop. This suggests the need to improve the parameter tuning approach further. }
	
	\begin{table}[h]
		\centering
		\caption{Average total scores of different techniques used for reoptimization.}
		\label{tab:finalscoresav}
		\begin{tabular}{p{0.5\textwidth}
				>{\raggedleft\arraybackslash}p{0.15\textwidth}
				>{\raggedleft\arraybackslash}p{0.15\textwidth}}
			\hline\noalign{\smallskip}
			Approach & Average total score & Improvement (\%)\\
			\noalign{\smallskip}\hline\noalign{\smallskip}
			Base & 0.727 & 0.00 \\
			Primal solution hint & 0.717 & 1.30 \\
			Reusing branching history & 0.712 & 2.09\\
			Reusing branching history with strong branching  & 0.708 & 2.57 \\
			Primal solution hint + Reusing branching history with strong branching & 0.695 & 4.40 \\
			(All above +) parameter tuning & 0.662 & 8.93\\
			(All above +) Other optimizations & 0.634 & 12.77\\
			\noalign{\smallskip}\hline
		\end{tabular}
	\end{table}
	
	\begin{table}[h]
		\centering
		\caption{Average total scores of reoptimization techniques per batch.}
		\label{tab:finalscoresbatchav}
		\begin{tabular}{lrrr}
			\hline\noalign{\smallskip}
			Batch & \multicolumn{2}{c}{Average total score} & Improvement (\%) \\
			\noalign{\smallskip}\cline{2-3}\noalign{\smallskip}
			& Base & Reoptimization &  \\
			\noalign{\smallskip}\hline\noalign{\smallskip}
			1-10 & 0.688 & 0.715 & -3.91 \\
			11-20 & 0.737 & 0.672 & 8.74 \\
			21-30 & 0.718 & 0.601 & 16.29 \\
			31-40 & 0.732 & 0.593 & 18.88 \\
			41-50 & 0.759 & 0.588 & 22.57 \\
			\noalign{\smallskip}\hline
		\end{tabular}
		
	\end{table}

	\begin{table}[h]
		\centering
		\caption{{Shifted geometric means of solving time for instances for each series.}}
		\label{tab:stgeometricmeans}
		\begin{tabular}{lrrr}
			\hline\noalign{\smallskip}
			Series & \multicolumn{2}{c}{Solve time geometric means} & Improvement (\%) \\
			\noalign{\smallskip}\cline{2-3}\noalign{\smallskip}
			& Base & Reoptimization &  \\
			\noalign{\smallskip}\hline\noalign{\smallskip}
			Bound 1 & 400.48 & 215.99 & 46.07 \\
			Bound 2 & 186.23 & 117.68 & 36.81 \\
			Bound 3 & 336.96 & 239.61 & 28.89 \\
			RHS 1 & 294.39 & 215.97 & 26.64 \\
			RHS 2 & 69.55 & 33.94 & 51.20 \\
			RHS 3 & 307.66 & 233.51 & 24.10 \\
			RHS 4 & 70 & 34.46 & 50.77 \\
			Objective 1 & 301.95 & 232.44 & 23.02 \\
			Objective 2 & 256.69 & 146.41 & 42.96 \\
			Objective 3 & 153.85 & 97.77 & 36.45 \\
			RHS Objective 1 & 492.8 & 319.52 & 35.16 \\
			RHS Objective 2 & 91.03 & 109 & -19.74 \\
			Matrix 1 & 270.7 & 298.94 & -10.43 \\
			Matrix RHS Objective Bound & 92.8 & 56.99 & 38.59 \\
			Matrix RHS Bound & 204.31 & 236.29 & -15.65\\
			\noalign{\smallskip}\hline
		\end{tabular}
		
	\end{table}
	
	\section{Conclusion and future work}
	\label{sec:conclusion}
	
	\begin{comment}
		\begin{itemize}
			\item Extension of primal solution.
			\item Extension of reusing pseudocosts (as branching priority)
			\item Extension of strong branching (increased reliability branching efforts)
			\item Tune more parameters.
			\item Using application specific knowledge to turn off more components.
		\end{itemize}
	\end{comment}
	
	In this paper, we described some techniques to strengthen and tune the MIP solver by reusing information generated by solving previous instances. The datasets we used to measure the impact contained 50 similar instances for each series. {These datasets were used for the MIP computational competition 2023}. This section summarizes our ideas and describes some ways to extend them.
	
	{In this paper, we tested our approach on a relatively small benchmark. In practice, the number of instances can be much larger than 50 for a given series. Moreover, all the instances in our benchmark can be solved within 10 minutes by any commercial solver, hence, they are relatively easy to solve. The competition committee introduced these limitations intentionally to limit the time required to experiment and evaluate various reoptimization techniques. Below, we suggest some ways to extend our approach to address these limitations.}
	
	First, we described how to reuse parts of previous instance solutions. Constructing a common partial solution and providing it as a solution hint works best in terms of converting the hints to a feasible solution. We use variable-value pairs that appear in at least $\alpha$\% of all the solutions to generate the common solution. One potential improvement to our approach is to dynamically change the value of $\alpha$ based on the hint-to-solution conversion numbers. This can potentially improve the conversion numbers. 
	
	As discussed in Section \ref{sec:finalcomparison}, most of the performance improvement comes from reusing the information to improve the dual bounds faster. Reusing the branching history alone brings more significant improvements than reusing the primal solutions. Pseudocosts can also be reused in the form of branching priorities for variables provided to the solver along with the model. We did not explore this direction.
	
	We use strong branching to solve the first instance to generate more reliable pseudocosts for the subsequent instances. {This technique resulted in significant gains in the RHS series. However, the effects of using strong branching in the first instance faded for the instances in later batches. In many series, using strong branching prevented us from proving optimality in the first instance.} An alternative approach can be to consider using reliability branching with increased thresholds for marking the pseudocost of a variable reliable. {If the time limits are not strict, we can also consider using more time for the first few instances to accommodate more strong branching than the default parameters.} The idea is to perform more strong branching in the initial instances to generate more reliable pseudocosts. Furthermore, we only use strong branching for the first instance. One can consider using more strong branching after solving every $n$ instances. 
	
	Automated parameter tuning is the most extensible technique that also leads to the most significant improvements in performance. Since we have only 50 instances in each series, our tuning approach focuses on faster convergence. We only tune 3 parameters as tuning more parameters simultaneously is not effective with our approach. However, if we have more instances, we can consider tuning more parameters sequentially, and more exploration (leading to slower but more reliable convergence) for each parameter tuned.
	
	Finally, we collect statistics of each major solver component and automatically disable some ineffective ones. For some series, this significantly improves the performance. However, the improvements are not consistent. 
	
	It is important to note that not {every} technique we discussed in this paper works well on each series. While each technique improves the average scores, they may significantly degrade the performance on some series. {Such mixed results can be partly explained by the complexity of the MIP solvers. MIP solvers have many components, and the interplay between them makes it hard to predict whether a particular technique will work well.}
	
	Except for reusing the branching history, the techniques discussed in this paper do not require access to the solver code {as m}ost solvers have public APIs for providing hints, statistics collections, and on/off control of various parts through parameters. We did not use any domain-specific knowledge in our approach apart from the limited use of the knowledge of the changing parts in each series.  {This means that the techniques we presented can be used for any reoptimization problem. Furthermore, we can get better results by using domain-specific knowledge along with these techniques.}

	%\begin{acknowledgements}
	%\end{acknowledgements}



	\bibliographystyle{spmpsci}
	\bibliography{references.bib}

\begin{thebibliography}{10}
\providecommand{\url}[1]{{#1}}
\providecommand{\urlprefix}{URL }
\expandafter\ifx\csname urlstyle\endcsname\relax
  \providecommand{\doi}[1]{DOI~\discretionary{}{}{}#1}\else
  \providecommand{\doi}{DOI~\discretionary{}{}{}\begingroup
  \urlstyle{rm}\Url}\fi

\bibitem{achterberg2005reliability}
Achterberg, T., Koch, T., Martin, A.: Branching rules revisited.
\newblock Operations Research Letters \textbf{33}(1), 42--54 (2005)

\bibitem{andreassian2021camels}
Andr{\'e}assian, V., Delaigue, O., Perrin, C., Janet, B., Addor, N.: Camels-fr:
  A large sample, hydroclimatic dataset for france, to support model testing
  and evaluation.
\newblock In: EGU General Assembly Conference Abstracts, pp. EGU21--13349
  (2021)

\bibitem{angulo2016improving}
Angulo, G., Ahmed, S., Dey, S.S.: Improving the integer l-shaped method.
\newblock INFORMS Journal on Computing \textbf{28}(3), 483--499 (2016)

\bibitem{applegate1995strongbranching}
Applegate, D., Bixby, R., Chv{\'a}tal, V., Cook, W.: Finding cuts in the TSP (A
  preliminary report), vol.~95.
\newblock Citeseer (1995)

\bibitem{auer2002finite}
Auer, P., Cesa-Bianchi, N., Fischer, P.: Finite-time analysis of the multiarmed
  bandit problem.
\newblock Machine learning \textbf{47}, 235--256 (2002)

\bibitem{ausiello2011complexity}
Ausiello, G., Bonifaci, V., Escoffier, B.: Complexity and approximation in
  reoptimization.
\newblock In: Computability in Context: Computation and Logic in the Real
  World, pp. 101--129. World Scientific (2011)

\bibitem{benichou1971pscost}
B{\'e}nichou, M., Gauthier, J.M., Girodet, P., Hentges, G., Ribi{\`e}re, G.,
  Vincent, O.: Experiments in mixed-integer linear programming.
\newblock Mathematical Programming \textbf{1}, 76--94 (1971)

\bibitem{berthold2006crossover}
Berthold, T.: Primal heuristics for mixed integer programs.
\newblock Ph.D. thesis, Zuse Institute Berlin (ZIB) (2006)

\bibitem{berthold2022restarttransfer}
Berthold, T., Hendel, G., Salvagnin, D.: Transferring information across
  restarts in mip.
\newblock In: Integration of Constraint Programming, Artificial Intelligence,
  and Operations Research: 19th International Conference, CPAIOR 2022, Los
  Angeles, CA, USA, June 20-23, 2022, Proceedings, pp. 24--33. Springer (2022)

\bibitem{BestuzhevaEtal2021ZR}
Bestuzheva, K., Besan{\c{c}}on, M., Chen, W.K., Chmiela, A., Donkiewicz, T.,
  van Doornmalen, J., Eifler, L., Gaul, O., Gamrath, G., Gleixner, A.,
  Gottwald, L., Graczyk, C., Halbig, K., Hoen, A., Hojny, C., van~der Hulst,
  R., Koch, T., L{\"u}bbecke, M., Maher, S.J., Matter, F., M{\"u}hmer, E.,
  M{\"u}ller, B., Pfetsch, M.E., Rehfeldt, D., Schlein, S., Schl{\"o}sser, F.,
  Serrano, F., Shinano, Y., Sofranac, B., Turner, M., Vigerske, S.,
  Wegscheider, F., Wellner, P., Weninger, D., Witzig, J.: {The SCIP
  Optimization Suite 8.0}.
\newblock ZIB-Report 21-41, Zuse Institute Berlin (2021).
\newblock \urlprefix\url{http://nbn-resolving.de/urn:nbn:de:0297-zib-85309}

\bibitem{blum2016cmsa}
Blum, C., Pinacho, P., L{\'o}pez-Ib{\'a}{\~n}ez, M., Lozano, J.A.: Construct,
  merge, solve \& adapt a new general algorithm for combinatorial optimization.
\newblock Computers \& Operations Research \textbf{68}, 75--88 (2016)

\bibitem{bolusani2023mip}
Bolusani, S., Besan{\c{c}}on, M., Gleixner, A., Berthold, T., d'Ambrosio, C.,
  Mu{\~n}oz, G., Paat, J., Thomopulos, D.: The mip workshop 2023 computational
  competition on reoptimization.
\newblock arXiv preprint arXiv:2311.14834  (2023)

\bibitem{bolusani2022framework}
Bolusani, S., Ralphs, T.K.: A framework for generalized benders’
  decomposition and its application to multilevel optimization.
\newblock Mathematical Programming \textbf{196}(1-2), 389--426 (2022)

\bibitem{sapse}
or~an SAP~affiliate company, S.S.: Milp benchmarks cellphoneco (2023).
\newblock
  \urlprefix\url{https://github.com/SAP-samples/ibp-sop-benchmarks-milp-cellphoneco}

\bibitem{cplex2020}
{Cplex, IBM ILOG}: {IBM ILOG CPLEX Optimizer} (2023).
\newblock
  \urlprefix\url{https://www.ibm.com/products/ilog-cplex-optimization-studio/cplex-optimizer}

\bibitem{danna2005rins}
Danna, E., Rothberg, E., Pape, C.L.: Exploring relaxation induced neighborhoods
  to improve {MIP} solutions.
\newblock Mathematical Programming \textbf{102}, 71--90 (2005)

\bibitem{firat2020column}
Firat, M., Crognier, G., Gabor, A.F., Hurkens, C.A., Zhang, Y.: Column
  generation based heuristic for learning classification trees.
\newblock Computers \& Operations Research \textbf{116}, 104866 (2020)

\bibitem{gamrath2019structure}
Gamrath, G., Berthold, T., Heinz, S., Winkler, M.: Structure-driven
  fix-and-propagate heuristics for mixed integer programming.
\newblock Mathematical Programming Computation \textbf{11}(4), 675--702 (2019)

\bibitem{gamrath2015reoptimization}
Gamrath, G., Hiller, B., Witzig, J.: Reoptimization techniques for {MIP}
  solvers.
\newblock In: Experimental Algorithms: 14th International Symposium, SEA 2015,
  Paris, France, June 29--July 1, 2015, Proceedings 14, pp. 181--192. Springer
  (2015)

\bibitem{gleixner2021miplib}
Gleixner, A., Hendel, G., Gamrath, G., Achterberg, T., Bastubbe, M., Berthold,
  T., Christophel, P., Jarck, K., Koch, T., Linderoth, J., et~al.: Miplib 2017:
  data-driven compilation of the 6th mixed-integer programming library.
\newblock Mathematical Programming Computation \textbf{13}(3), 443--490 (2021)

\bibitem{gurobi}
{Gurobi Optimization, LLC}: {Gurobi Optimizer Reference Manual} (2023).
\newblock \urlprefix\url{https://www.gurobi.com}

\bibitem{hendel2015pscostvariance}
Hendel, G.: Enhancing {MIP} branching decisions by using the sample variance of
  pseudo costs.
\newblock In: Integration of AI and OR Techniques in Constraint Programming:
  12th International Conference, CPAIOR 2015, Barcelona, Spain, May 18-22,
  2015, Proceedings 12, pp. 199--214. Springer (2015)

\bibitem{hendel2022adaptive}
Hendel, G.: Adaptive large neighborhood search for mixed integer programming.
\newblock Mathematical Programming Computation pp. 1--37 (2022)

\bibitem{jimenez2022warm}
Jim{\'e}nez-Cordero, A., Morales, J.M., Pineda, S.: Warm-starting constraint
  generation for mixed-integer optimization: A machine learning approach.
\newblock Knowledge-Based Systems \textbf{253}, 109570 (2022)

\bibitem{li2017hyperband}
Li, L., Jamieson, K., DeSalvo, G., Rostamizadeh, A., Talwalkar, A.: Hyperband:
  A novel bandit-based approach to hyperparameter optimization.
\newblock The journal of machine learning research \textbf{18}(1), 6765--6816
  (2017)

\bibitem{linderoth1999searchstrategies}
Linderoth, J.T., Savelsbergh, M.W.: A computational study of search strategies
  for mixed integer programming.
\newblock INFORMS Journal on Computing \textbf{11}(2), 173--187 (1999)

\bibitem{ralphs2006duality}
Ralphs, T., G{\"u}zelsoy, M.: Duality and warm starting in integer programming.
\newblock In: The proceedings of the 2006 NSF design, service, and
  manufacturing grantees and research conference (2006)

\bibitem{sabharwal2012guiding}
Sabharwal, A., Samulowitz, H., Reddy, C.: Guiding combinatorial optimization
  with uct.
\newblock In: Integration of AI and OR Techniques in Contraint Programming for
  Combinatorial Optimzation Problems: 9th International Conference, CPAIOR
  2012, Nantes, France, May 28--June1, 2012. Proceedings 9, pp. 356--361.
  Springer (2012)

\bibitem{schieber2018theory}
Schieber, B., Shachnai, H., Tamir, G., Tamir, T.: A theory and algorithms for
  combinatorial reoptimization.
\newblock Algorithmica \textbf{80}, 576--607 (2018)

\bibitem{tanner2010iis}
Tanner, M.W., Ntaimo, L.: Iis branch-and-cut for joint chance-constrained
  stochastic programs and application to optimal vaccine allocation.
\newblock European Journal of Operational Research \textbf{207}(1), 290--296
  (2010)

\bibitem{thomopulos2023generating}
Thomopulos, D., van Ackooij, W., D’Ambrosio, C., St{\'e}fanon, M.: Generating
  hydro unit commitment instances.
\newblock TOP pp. 1--31 (2023)

\bibitem{witzig2014reoptimization}
Witzig, J.: Reoptimization techniques in {MIP} solvers, master's thesis, {TU
  Berlin} (2014)

\bibitem{xu2011hydra}
Xu, L., Hutter, F., Hoos, H.H., Leyton-Brown, K.: Hydra-mip: Automated
  algorithm configuration and selection for mixed integer programming.
\newblock In: RCRA workshop on experimental evaluation of algorithms for
  solving problems with combinatorial explosion at the international joint
  conference on artificial intelligence (IJCAI), pp. 16--30 (2011)

\end{thebibliography}
	
	\section*{Statements and Declarations}
	
	\subsection*{Funding}
	
	Partial financial support was received from Canada Excellence Research Chair (CERC), Polytechnique Montreal. The computational resources
	for experimental evaluation were provided by CERC.
	%The authors declare that no funds, grants, or other support were received during the preparation of this manuscript.
	
	\subsection*{Competing Interests}
	
	The authors have no relevant financial or non-financial interests to disclose.
	
	\subsection*{Author Contributions}
	
	Krunal Kishor Patel designed, implemented, and tested the ideas presented in this paper.
	
	\subsection*{Data Availability}
	
	The datasets analysed during the current study are available in the github repository, https://github.com/krooonal/mipcc2023/tree/MPCPaper
	
	\appendix
	
	\section{Detailed Results}

	\begin{longtable}{l
			>{\raggedleft\arraybackslash}p{0.12\textwidth}
			>{\raggedleft\arraybackslash}p{0.08\textwidth}
			>{\raggedleft\arraybackslash}p{0.1\textwidth}
			>{\raggedleft\arraybackslash}p{0.12\textwidth}}
		\caption{Effects of providing solution hints on total scores. {Best scores are in bold.}}
		\label{tab:solutionhintscores}\\
		\hline\noalign{\smallskip}
		Series & Conversion rate (\%) & Batch & \multicolumn{2}{c}{Total score} \\
		\noalign{\smallskip}\cline{4-5}\noalign{\smallskip}
		&  &  & {BASE} & {SOLHINT} \\
		\noalign{\smallskip}\hline\noalign{\smallskip}
		\multirow{5}{*}{Bound 1} & \multirow{5}{*}{51.0} & 1-10 & \textbf{0.78} & 0.91 \\
		&  & 11-20 & \textbf{0.77} & 0.83 \\
		&  & 21-30 & \textbf{0.72} & 0.78 \\
		&  & 31-40 & \textbf{0.80} & 0.85 \\
		&  & 41-50 & \textbf{0.60} & 0.70 \\
		\noalign{\smallskip}\hline\noalign{\smallskip}
		\multirow{5}{*}{Bound 2} & \multirow{5}{*}{42.9} & 1-10 & \textbf{0.76} & 0.78 \\
		&  & 11-20 & 0.60 & \textbf{0.58} \\
		&  & 21-30 & \textbf{0.64} & 0.66 \\
		&  & 31-40 & 0.66 & \textbf{0.59} \\
		&  & 41-50 & \textbf{0.63} & 0.69 \\
		\noalign{\smallskip}\hline\noalign{\smallskip}
		\multirow{5}{*}{Bound 3} & \multirow{5}{*}{61.2} & 1-10 & 0.62 & \textbf{0.61} \\
		&  & 11-20 & 0.68 & \textbf{0.67} \\
		&  & 21-30 & 0.70 & \textbf{0.69} \\
		&  & 31-40 & 0.61 & \textbf{0.60} \\
		&  & 41-50 & 0.57 & \textbf{0.53} \\
		\noalign{\smallskip}\hline\noalign{\smallskip}
		\multirow{5}{*}{RHS 1} & \multirow{5}{*}{69.4} & 1-10 & \textbf{0.53} & 0.56 \\
		&  & 11-20 & 0.69 & \textbf{0.67} \\
		&  & 21-30 & 0.72 & \textbf{0.71} \\
		&  & 31-40 & 0.89 & \textbf{0.79} \\
		&  & 41-50 & \textbf{0.90} & 0.91 \\
		\noalign{\smallskip}\hline\noalign{\smallskip}
		\multirow{5}{*}{RHS 2} & \multirow{5}{*}{100.0} & 1-10 & 0.99 & \textbf{0.96} \\
		&  & 11-20 & 0.99 & {0.99} \\
		&  & 21-30 & 1.00 & {1.00} \\
		&  & 31-40 & 1.00 & \textbf{0.98} \\
		&  & 41-50 & 0.99 & \textbf{0.90} \\
		\noalign{\smallskip}\hline\noalign{\smallskip}
		\multirow{5}{*}{RHS 3} & \multirow{5}{*}{0.0} & 1-10 & {0.66} & 0.66 \\
		&  & 11-20 & \textbf{0.78} & 0.79 \\
		&  & 21-30 & {0.71} & 0.71 \\
		&  & 31-40 & {0.54} & 0.54 \\
		&  & 41-50 & {0.57} & 0.57 \\
		\noalign{\smallskip}\hline\noalign{\smallskip}
		\multirow{5}{*}{RHS 4} & \multirow{5}{*}{100.0} & 1-10 & 1.00 & \textbf{0.96} \\
		&  & 11-20 & 1.00 & {1.00} \\
		&  & 21-30 & 1.00 & \textbf{0.99} \\
		&  & 31-40 & 1.00 & \textbf{0.98} \\
		&  & 41-50 & 1.00 & \textbf{0.99} \\
		\noalign{\smallskip}\hline\noalign{\smallskip}
		\multirow{5}{*}{Objective 1} & \multirow{5}{*}{100.0} & 1-10 & 0.49 & \textbf{0.48} \\
		&  & 11-20 & 0.59 & \textbf{0.55} \\
		&  & 21-30 & 0.77 & \textbf{0.75} \\
		&  & 31-40 & 0.93 & \textbf{0.88} \\
		&  & 41-50 & 1.00 & \textbf{0.95} \\
		\noalign{\smallskip}\hline\noalign{\smallskip}
		\multirow{5}{*}{Objective 2} & \multirow{5}{*}{100.0} & 1-10 & 0.82 & {0.82} \\
		&  & 11-20 & \textbf{0.97} & 1.03 \\
		&  & 21-30 & \textbf{0.97} & 1.04 \\
		&  & 31-40 & 0.86 & \textbf{0.84} \\
		&  & 41-50 & \textbf{0.93} & 0.98 \\
		\noalign{\smallskip}\hline\noalign{\smallskip}
		\multirow{5}{*}{Objective 3} & \multirow{5}{*}{100.0} & 1-10 & {0.36} & 0.36 \\
		&  & 11-20 & 0.37 & {0.37} \\
		&  & 21-30 & \textbf{0.46} & 0.57 \\
		&  & 31-40 & 0.44 & \textbf{0.38} \\
		&  & 41-50 & 0.47 & \textbf{0.45} \\
		\noalign{\smallskip}\hline\noalign{\smallskip}
		\multirow{5}{*}{RHS Objective 1} & \multirow{5}{*}{93.9} & 1-10 & 1.01 & \textbf{0.91} \\
		&  & 11-20 & \textbf{0.99} & 1.00 \\
		&  & 21-30 & \textbf{0.97} & 0.99 \\
		&  & 31-40 & 0.93 & \textbf{0.92} \\
		&  & 41-50 & 1.08 & \textbf{0.94} \\
		\noalign{\smallskip}\hline\noalign{\smallskip}
		\multirow{5}{*}{RHS Objective 2} & \multirow{5}{*}{89.8} & 1-10 & 0.41 & \textbf{0.33} \\
		&  & 11-20 & \textbf{0.46} & 0.50 \\
		&  & 21-30 & \textbf{0.23} & 0.32 \\
		&  & 31-40 & 0.40 & \textbf{0.34} \\
		&  & 41-50 & {0.59} & 0.59 \\
		\noalign{\smallskip}\hline\noalign{\smallskip}
		\multirow{5}{*}{Matrix 1} & \multirow{5}{*}{100.0} & 1-10 & \textbf{0.82} & 0.91 \\
		&  & 11-20 & 0.94 & \textbf{0.92} \\
		&  & 21-30 & \textbf{0.90} & 0.95 \\
		&  & 31-40 & \textbf{0.89} & 0.98 \\
		&  & 41-50 & \textbf{0.93} & 1.01 \\
		\noalign{\smallskip}\hline\noalign{\smallskip}
		\multirow{5}{*}{Matrix RHS Objective Bound} & \multirow{5}{*}{87.8} & 1-10 & 0.60 & \textbf{0.44} \\
		&  & 11-20 & 0.55 & \textbf{0.38} \\
		&  & 21-30 & 0.50 & \textbf{0.30} \\
		&  & 31-40 & 0.36 & \textbf{0.33} \\
		&  & 41-50 & 0.51 & \textbf{0.36} \\
		\noalign{\smallskip}\hline\noalign{\smallskip}
		\multirow{5}{*}{Matrix RHS Bound} & \multirow{5}{*}{0.0} & 1-10 & 0.47 & {0.47} \\
		&  & 11-20 & 0.66 & {0.66} \\
		&  & 21-30 & 0.48 & {0.48} \\
		&  & 31-40 & 0.66 & {0.66} \\
		&  & 41-50 & 0.59 & {0.59} \\
		\noalign{\smallskip}\hline
	\end{longtable}
	
	\begin{longtable}{l
			>{\raggedleft\arraybackslash}p{0.08\textwidth}
			>{\raggedleft\arraybackslash}p{0.1\textwidth}
			>{\raggedleft\arraybackslash}p{0.1\textwidth}
			>{\raggedleft\arraybackslash}p{0.12\textwidth}}
		\caption{Effects of reusing branching history and strong branching for the first instance on total scores. {Best scores are in bold.}}
		\label{tab:reusehistoryscores}\\
		\hline\noalign{\smallskip}
		Series &Batch &  \multicolumn{3}{c}{Total score}  \\
		\noalign{\smallskip}\cline{3-5}\noalign{\smallskip}
		& & BASE & REHIS & REHIS+SB \\
		\noalign{\smallskip}\hline\noalign{\smallskip}
		\multirow{5}{*}{Bound 1} &1-10& 0.78 & \textbf{0.62} & 0.64 \\
		&11-20& 0.77 & \textbf{0.46} & \textbf{0.46} \\
		&21-30& 0.72 & \textbf{0.49} & \textbf{0.49} \\
		&31-40& 0.80 & \textbf{0.51} & 0.57 \\
		&41-50& 0.60 & \textbf{0.34} & 0.40 \\
		\noalign{\smallskip}\hline\noalign{\smallskip}
		\multirow{5}{*}{Bound 2}&1-10 & 0.76 & \textbf{0.61} & \textbf{0.61} \\
		&11-20& \textbf{0.60} & {0.61} & 0.63 \\
		&21-30& 0.64 & \textbf{0.53} & 0.62 \\
		&31-40& 0.66 & 0.60 & \textbf{0.55} \\
		&41-50& 0.63 & \textbf{0.52} & 0.68 \\
		\noalign{\smallskip}\hline\noalign{\smallskip}
		\multirow{5}{*}{Bound 3}&1-10 & 0.62 & \textbf{0.61} & 0.66 \\
		&11-20& 0.68 & \textbf{0.67} & \textbf{0.67} \\
		&21-30& 0.70 & 0.72 & \textbf{0.68} \\
		&31-40& 0.61 & \textbf{0.54} & 0.63 \\
		&41-50& 0.57 & 0.49 & \textbf{0.44} \\
		\noalign{\smallskip}\hline\noalign{\smallskip}
		\multirow{5}{*}{RHS 1} &1-10& 0.53 & \textbf{0.51} & 0.62 \\
		&11-20& 0.69 & \textbf{0.64} & 0.70 \\
		&21-30& 0.72 & 0.77 & \textbf{0.62} \\
		&31-40& 0.89 & 0.77 & \textbf{0.71} \\
		&41-50& 0.90 & 0.80 & \textbf{0.79} \\
		\noalign{\smallskip}\hline\noalign{\smallskip}
		\multirow{5}{*}{RHS 2}&1-10 & 0.99 & 0.99 &\textbf{0.97} \\
		&11-20& 0.99 & 0.99 & \textbf{0.92} \\
		&21-30& 1.00 & 1.00 & \textbf{0.97} \\
		&31-40& 1.00 & 0.99 & \textbf{0.93} \\
		&41-50& 0.99 & 1.00 & \textbf{0.85} \\
		\noalign{\smallskip}\hline\noalign{\smallskip}
		\multirow{5}{*}{RHS 3}&1-10 & 0.66 & 0.62 & \textbf{0.56} \\
		&11-20& 0.78 & \textbf{0.51} & 0.57 \\
		&21-30& 0.71 & 0.69 & \textbf{0.62} \\
		&31-40& 0.54 & 0.59 & \textbf{0.52} \\
		&41-50& 0.57 & 0.59 & \textbf{0.52} \\
		\noalign{\smallskip}\hline\noalign{\smallskip}
		\multirow{5}{*}{RHS 4}&1-10 & 1.00 & 1.00 & \textbf{0.94} \\
		&11-20& 1.00 & 1.00 & \textbf{0.94} \\
		&21-30& 1.00 & 1.00 & \textbf{0.96} \\
		&31-40& 1.00 & 1.00 & \textbf{0.95} \\
		&41-50& 1.00 & 1.00 & \textbf{0.92} \\
		\noalign{\smallskip}\hline\noalign{\smallskip}
		\multirow{5}{*}{Objective 1}&1-10 & \textbf{0.49} & \textbf{0.49} & 0.52 \\
		&11-20& 0.59 & \textbf{0.57} & \textbf{0.57} \\
		&21-30& \textbf{0.77} & \textbf{0.77} & 0.78 \\
		&31-40& 0.93 & 0.89 & \textbf{0.87} \\
		&41-50& 1.00 & 0.97 & \textbf{0.92} \\
		\noalign{\smallskip}\hline\noalign{\smallskip}
		\multirow{5}{*}{Objective 2}&1-10 & \textbf{0.82} & 0.92 & 1.08 \\
		&11-20& \textbf{0.97} & 1.05 & 1.05 \\
		&21-30& 0.97 & \textbf{0.95} & 0.96 \\
		&31-40& 0.86 & \textbf{0.81} & 0.85 \\
		&41-50& \textbf{0.93} & 1.10 & 1.12 \\
		\noalign{\smallskip}\hline\noalign{\smallskip}
		\multirow{5}{*}{Objective 3}&1-10 & \textbf{0.36} & 0.38 & 0.59 \\
		&11-20& \textbf{0.37} & 0.57 & 0.54 \\
		&21-30& 0.46 & 0.56 & \textbf{0.44} \\
		&31-40& \textbf{0.44} & 0.48 & \textbf{0.44} \\
		&41-50& 0.47 & \textbf{0.37} & 0.42 \\
		\noalign{\smallskip}\hline\noalign{\smallskip}
		\multirow{5}{*}{RHS Objective 1}&1-10 & 1.01 & \textbf{0.96} & 1.08 \\
		&11-20& \textbf{0.99} & 1.02 & 1.04 \\
		&21-30& \textbf{0.97} & 1.00 & 1.01 \\
		&31-40& \textbf{0.93} & 0.99 & \textbf{0.93} \\
		&41-50& 1.08 & \textbf{1.00} & \textbf{1.00} \\
		\noalign{\smallskip}\hline\noalign{\smallskip}
		\multirow{5}{*}{RHS Objective 2}&1-10 & \textbf{0.41} & 0.42 & 0.47 \\
		&11-20& \textbf{0.46} & 0.48 & 0.49 \\
		&21-30& \textbf{0.23} & 0.37 & {0.24} \\
		&31-40& \textbf{0.40} & 0.44 & \textbf{0.40} \\
		&41-50& 0.59 & 0.58 & \textbf{0.56} \\
		\noalign{\smallskip}\hline\noalign{\smallskip}
		\multirow{5}{*}{Matrix 1}&1-10 & \textbf{0.81} & 0.87 & 0.90 \\
		&11-20& 0.94 & \textbf{0.85} & 0.97 \\
		&21-30& \textbf{0.90} & 0.96 & 0.98 \\
		&31-40& \textbf{0.89} & 0.98 & 0.94 \\
		&41-50& 0.93 & 0.89 & \textbf{0.87} \\
		\noalign{\smallskip}\hline\noalign{\smallskip}
		\multirow{5}{*}{Matrix RHS Objective Bound}&1-10 & \textbf{0.60} & 0.68 & 0.75 \\
		&11-20& \textbf{0.55} & 0.63 & 0.60 \\
		&21-30& \textbf{0.50} & 0.65 & 0.63 \\
		&31-40& \textbf{0.36} & 0.46 & 0.43 \\
		&41-50& \textbf{0.51} & 0.59 & 0.54 \\
		\noalign{\smallskip}\hline\noalign{\smallskip}
		\multirow{5}{*}{Matrix RHS Bound}&1-10 & \textbf{0.47} & 0.59 & 0.53 \\
		&11-20& 0.67 & \textbf{0.54} & 0.71 \\
		&21-30& 0.48 & 0.57 & \textbf{0.45} \\
		&31-40& 0.66 & \textbf{0.53} & 0.60 \\
		&41-50& 0.60 & 0.67 & \textbf{0.55} \\
		\noalign{\smallskip}\hline
	\end{longtable}
	
		\begin{longtable}{l
				>{\raggedleft\arraybackslash}p{0.15\textwidth}
				>{\raggedleft\arraybackslash}p{0.17\textwidth}
				>{\raggedleft\arraybackslash}p{0.16\textwidth}}
			\caption{Effects of automated parameter tuning on total scores. {Best scores are in bold.}}
			\label{tab:paramtuningscores}\\
			\hline\noalign{\smallskip}
			Series & Batch & \multicolumn{2}{c}{Total scores} \\
			\noalign{\smallskip}\cline{3-4}\noalign{\smallskip}
			& & NOTUNING & TUNING \\
			\noalign{\smallskip}\hline\noalign{\smallskip}
			\multirow{5}{*}{Bound 1} &1-10& \textbf{0.73} & 0.85 \\
			&11-20& 0.61 & \textbf{0.50} \\
			&21-30& 0.71 & \textbf{0.45} \\
			&31-40& 0.84 & \textbf{0.44} \\
			&41-50& 0.61 & \textbf{0.32} \\
			\noalign{\smallskip}\hline\noalign{\smallskip}
			\multirow{5}{*}{Bound 2}&1-10 & \textbf{0.57} & 0.58 \\
			&11-20& \textbf{0.45} & 0.50 \\
			&21-30& \textbf{0.40} & 0.43 \\
			&31-40& 0.40 & \textbf{0.39} \\
			&41-50& 0.40 & {0.40} \\
			\noalign{\smallskip}\hline\noalign{\smallskip}
			\multirow{5}{*}{Bound 3}&1-10 & \textbf{0.54} & 0.58 \\
			&11-20& \textbf{0.48} & 0.62 \\
			&21-30& \textbf{0.52} & 0.64 \\
			&31-40& \textbf{0.49} & 0.50 \\
			&41-50& \textbf{0.39} & 0.42 \\
			\noalign{\smallskip}\hline\noalign{\smallskip}
			\multirow{5}{*}{RHS 1}&1-10 & \textbf{0.50} & 0.77 \\
			&11-20& \textbf{0.64} & 0.68 \\
			&21-30& 0.71 & \textbf{0.61} \\
			&31-40& 0.72 & \textbf{0.69} \\
			&41-50& 0.79 & \textbf{0.48} \\
			\noalign{\smallskip}\hline\noalign{\smallskip}
			\multirow{5}{*}{RHS 2}&1-10 & \textbf{0.66} & 0.68 \\
			&11-20& 0.63 & \textbf{0.55} \\
			&21-30& 0.70 & \textbf{0.57} \\
			&31-40& 0.67 & \textbf{0.59} \\
			&41-50& 0.66 & \textbf{0.52} \\
			\noalign{\smallskip}\hline\noalign{\smallskip}
			\multirow{5}{*}{RHS 3}&1-10 & \textbf{0.61} & 0.79 \\
			&11-20& \textbf{0.75} & 0.85 \\
			&21-30& 0.57 & \textbf{0.55} \\
			&31-40& \textbf{0.52} & 0.57 \\
			&41-50& 0.53 & \textbf{0.37} \\
			\noalign{\smallskip}\hline\noalign{\smallskip}
			\multirow{5}{*}{RHS 4}&1-10 & 0.71 & \textbf{0.69} \\
			&11-20& 0.64 & \textbf{0.58} \\
			&21-30& 0.70 & \textbf{0.54} \\
			&31-40& 0.67 & \textbf{0.56} \\
			&41-50& 0.67 & \textbf{0.56} \\
			\noalign{\smallskip}\hline\noalign{\smallskip}
			\multirow{5}{*}{Objective 1} &1-10& \textbf{0.51} & 0.52 \\
			&11-20& 0.56 & \textbf{0.49} \\
			&21-30& 0.78 & \textbf{0.60} \\
			&31-40& 0.89 & \textbf{0.74} \\
			&41-50& 0.98 & \textbf{0.96} \\
			\noalign{\smallskip}\hline\noalign{\smallskip}
			\multirow{5}{*}{Objective 2} &1-10& 0.97 & \textbf{0.93} \\
			&11-20& \textbf{1.09} & 1.20 \\
			&21-30& 1.11 & \textbf{1.07} \\
			&31-40& \textbf{0.80} & 0.91 \\
			&41-50& 1.07 & \textbf{1.04} \\
			\noalign{\smallskip}\hline\noalign{\smallskip}
			\multirow{5}{*}{Objective 3}&1-10 & 0.46 & \textbf{0.38} \\
			&11-20& 0.63 & \textbf{0.37} \\
			&21-30& 0.64 & \textbf{0.31} \\
			&31-40& 0.79 & \textbf{0.45} \\
			&41-50& 0.76 & \textbf{0.50} \\
			\noalign{\smallskip}\hline\noalign{\smallskip}
			\multirow{5}{*}{RHS Objective 1}&1-10 & 1.07 & \textbf{1.03} \\
			&11-20& 1.00 & \textbf{0.98} \\
			&21-30& 1.00 & \textbf{0.89} \\
			&31-40& 0.93 & 0.93 \\
			&41-50& 0.96 & \textbf{0.83} \\
			\noalign{\smallskip}\hline\noalign{\smallskip}
			\multirow{5}{*}{RHS Objective 2} &1-10& \textbf{0.49} & 0.77 \\
			&11-20& \textbf{0.59} & 0.77 \\
			&21-30& \textbf{0.53} & 0.62 \\
			&31-40& \textbf{0.59} & 0.69 \\
			&41-50& \textbf{0.66} & 0.80 \\
			\noalign{\smallskip}\hline\noalign{\smallskip}
			\multirow{5}{*}{Matrix 1} &1-10& \textbf{1.00} & 1.04 \\
			&11-20& \textbf{0.95} & 1.02 \\
			&21-30& \textbf{1.00} & 1.03 \\
			&31-40& \textbf{1.05} & 1.06 \\
			&41-50& 1.08 & \textbf{1.06} \\
			\noalign{\smallskip}\hline\noalign{\smallskip}
			\multirow{5}{*}{Matrix RHS Objective Bound}&1-10 & \textbf{0.63} & 0.66 \\
			&11-20& 0.48 & \textbf{0.44} \\
			&21-30& 0.30 & \textbf{0.27} \\
			&31-40& 0.42 & \textbf{0.33} \\
			&41-50& 0.50 & \textbf{0.30} \\
			\noalign{\smallskip}\hline\noalign{\smallskip}
			\multirow{5}{*}{Matrix RHS Bound}&1-10 & \textbf{0.61} & 0.91 \\
			&11-20& \textbf{0.64} & 0.75 \\
			&21-30& 0.75 & 0.75 \\
			&31-40& 0.90 & \textbf{0.77} \\
			&41-50& 0.73 & \textbf{0.69} \\
			\noalign{\smallskip}\hline
		\end{longtable}

		\begin{longtable}{l
				>{\raggedleft\arraybackslash}p{0.15\textwidth}
				>{\raggedleft\arraybackslash}p{0.20\textwidth}
				>{\raggedleft\arraybackslash}p{0.16\textwidth}}
		\caption{Effects of turning off inefficient solver parts on total scores. {Best scores are in bold.}}
		\label{tab:turnoffscores}\\
		\hline\noalign{\smallskip}
		Series & Batch & \multicolumn{2}{c}{Total scores} \\
		\noalign{\smallskip}\cline{3-4}\noalign{\smallskip}
		&  & NOTURNOFF & TURNOFF \\
		\noalign{\smallskip}\hline\noalign{\smallskip}
		\multirow{5}{*}{Bound 1} & 1-10 & 0.78 & 0.78 \\
		& 11-20 & 0.48 & \textbf{0.46} \\
		& 21-30 & 0.46 &\textbf{0.42} \\
		& 31-40 & 0.48 & \textbf{0.39} \\
		& 41-50 & \textbf{0.34} & 0.35 \\
		\noalign{\smallskip}\hline\noalign{\smallskip}
		\multirow{5}{*}{Bound 2} & 1-10 & 0.59 &  0.59 \\
		& 11-20 & \textbf{0.43} & 0.50 \\
		& 21-30 & \textbf{0.38} & 0.42 \\
		& 31-40 & 0.36 & 0.36 \\
		& 41-50 & 0.41 & 0.41 \\
		\noalign{\smallskip}\hline\noalign{\smallskip}
		\multirow{5}{*}{Bound 3} & 1-10 & 0.58 & 0.58 \\
		& 11-20 & 0.59 & 0.59 \\
		& 21-30 & \textbf{0.54} & 0.56 \\
		& 31-40 & \textbf{0.49} & 0.52 \\
		& 41-50 & \textbf{0.40} & 0.42 \\
		\noalign{\smallskip}\hline\noalign{\smallskip}
		\multirow{5}{*}{RHS 1} & 1-10 & 0.76 & 0.76 \\
		& 11-20 & 0.71 & \textbf{0.62} \\
		& 21-30 & \textbf{0.55} & 0.63 \\
		& 31-40 & \textbf{0.62} & 0.63 \\
		& 41-50 & 0.66 & \textbf{0.54} \\
		\noalign{\smallskip}\hline\noalign{\smallskip}
		\multirow{5}{*}{RHS 2} & 1-10 & 0.66 & 0.66 \\
		& 11-20 & 0.54 & 0.54 \\
		& 21-30 & 0.59 &\textbf{0.42} \\
		& 31-40 & 0.55 & \textbf{0.27} \\
		& 41-50 & 0.51 & \textbf{0.26} \\
		\noalign{\smallskip}\hline\noalign{\smallskip}
		\multirow{5}{*}{RHS 3} & 1-10 &  0.62 &  0.62 \\
		& 11-20 & \textbf{0.72} & 0.82 \\
		& 21-30 & \textbf{0.50} & 0.65 \\
		& 31-40 & \textbf{0.62} & 0.79 \\
		& 41-50 & \textbf{0.46} & 0.51 \\
		\noalign{\smallskip}\hline\noalign{\smallskip}
		\multirow{5}{*}{RHS 4} & 1-10 & 0.70 & 0.70 \\
		& 11-20 & 0.57 & 0.57 \\
		& 21-30 & 0.53 & \textbf{0.42} \\
		& 31-40 & 0.56 & \textbf{0.30} \\
		& 41-50 & 0.55 & \textbf{0.29} \\
		\noalign{\smallskip}\hline\noalign{\smallskip}
		\multirow{5}{*}{Objective 1} & 1-10 & 0.48 & 0.48 \\
		& 11-20 & 0.45 & \textbf{0.43} \\
		& 21-30 & 0.67 & \textbf{0.64} \\
		& 31-40 & 0.74 & \textbf{0.72} \\
		& 41-50 & 0.90 & \textbf{0.88} \\
		\noalign{\smallskip}\hline\noalign{\smallskip}
		\multirow{5}{*}{Objective 2} & 1-10 & 0.93 & 0.93 \\
		& 11-20 & 1.29 & \textbf{1.18} \\
		& 21-30 & 1.00 & \textbf{0.95} \\
		& 31-40 & \textbf{0.94} & 0.95 \\
		& 41-50 & \textbf{1.07} & 1.13 \\
		\noalign{\smallskip}\hline\noalign{\smallskip}
		\multirow{5}{*}{Objective 3} & 1-10 & 0.43 &  0.43 \\
		& 11-20 & \textbf{0.39} & 0.49 \\
		& 21-30 & \textbf{0.40} & 0.46 \\
		& 31-40 & 0.41 & \textbf{0.35} \\
		& 41-50 & \textbf{0.36} & 0.42 \\
		\noalign{\smallskip}\hline\noalign{\smallskip}
		\multirow{5}{*}{RHS Objective 1} & 1-10 & 1.03 &  1.03 \\
		& 11-20 & 1.03 & \textbf{1.00} \\
		& 21-30 & \textbf{0.89} & 0.95 \\
		& 31-40 & \textbf{0.72} & 0.86 \\
		& 41-50 & 0.82 & 0.82 \\
		\noalign{\smallskip}\hline\noalign{\smallskip}
		\multirow{5}{*}{RHS Objective 2}& 1-10 & 0.76 &  0.76 \\
		& 11-20 & \textbf{0.63} & 0.74 \\
		& 21-30 & 0.69 & \textbf{0.56} \\
		& 31-40 & 0.63 & \textbf{0.50} \\
		& 41-50 & 0.72 & \textbf{0.68} \\
		\noalign{\smallskip}\hline\noalign{\smallskip}
		\multirow{5}{*}{Matrix 1}& 1-10 & 1.04 & 1.04 \\
		& 11-20 & \textbf{0.97} & 1.04 \\
		& 21-30 & 1.05 & 1.05 \\
		& 31-40 & \textbf{1.06} & 1.07 \\
		& 41-50 & 1.05 & \textbf{1.04} \\
		\noalign{\smallskip}\hline\noalign{\smallskip}
		\multirow{5}{*}{Matrix RHS Objective Bound}  & 1-10 & 0.58 & 0.58 \\
		& 11-20 & 0.35 & 0.35 \\
		& 21-30 & \textbf{0.25} & 0.28 \\
		& 31-40 & \textbf{0.32} & 0.33 \\
		& 41-50 & 0.30 & 0.30 \\
		\noalign{\smallskip}\hline\noalign{\smallskip}
		\multirow{5}{*}{Matrix RHS Bound} & 1-10 & 0.80 & 0.80 \\
		& 11-20 & 0.78 & 0.78 \\
		& 21-30 & 0.67 & \textbf{0.61} \\
		& 31-40 & 0.93 & \textbf{0.86} \\
		& 41-50 & \textbf{0.65} & 0.76 \\
		\noalign{\smallskip}\hline
		%Average & & 0.662 & 0.634\\
		%\hline
	\end{longtable}

\end{document}